\documentclass[smallextended,natbib,runningheads]{svjour3}
\journalname{Annals of the Institute of Statistical Mathematics}
\smartqed
\usepackage{amsmath,amssymb}
\usepackage[left=3.8cm,right=3.8cm]{geometry}

\newcommand{\bigsum}{\displaystyle \sum}
\newcommand{\W}{W_{(i)}}
\newcommand{\round}{\lceil np \rceil}
\newcommand{\roundp}{\lceil n(1-p) \rceil}
\newcommand{\rround}{\lceil nr_n \rceil}
\newcommand{\less}{\mathcal{I}_{W_i<F^{-1}(p)}}
\newcommand{\more}{\mathcal{I}_{W_i\geq F^{-1}(p)}}
\newcommand{\lessone}{\mathcal{I}_{W_1<F^{-1}(p)}}
\newcommand{\moreone}{\mathcal{I}_{W_1\geq F^{-1}(p)}}

\begin{document}
\title{Asymptotics of the Empirical Cross-over Function}
\author{Karthik Bharath \and Vladimir Pozdnyakov \and Dipak. K. Dey
}
\institute{University of Connecticut\\
215 Glenbrook Road, Storrs, CT-06269}
\date{}
\maketitle
\begin{abstract}
We consider a combination of heavily trimmed sums and sample quantiles which arises when examining  properties of clustering criteria and prove limit theorems. The object of interest, which we call the Empirical Cross-over Function, is an L-statistic whose weights do not comply with the requisite regularity conditions for usage of existing limit results. The law of large numbers, CLT and a functional CLT are proven.
\end{abstract}
\keywords{Clustering \and L-statistics \and CLT \and Functional CLT}
\PACS{62F05 \and 62G30 \and 62E20}

\section{Introduction}

Suppose $W_1,W_2,\dots,W_n$ for $n \geq 1$  are i.i.d. random variables with distribution function $F$. If $W_{(1)}\leq W_{(2)} \leq \cdots \leq W_{(n)}$ are the order statistics, then, we define, for $0<p<1$, the \emph{Empirical Cross-over Function} (ECF)
\begin{equation}\label{G_function}
 G_n(p)=\frac{1}{k}\sum_{j=1}^k W_{(j)}-W_{(k)} +\frac{1}{n-k}\sum_{j=k+1}^n W_{(j)}-W_{(k+1)} \qquad \textrm{for  }\frac{k-1}{n}\leq p <\frac{k}{n}.
\end{equation}
The function $G_n$ is a special case of linear functions of order statistics $W_{(i)}$, $1 \leq i \leq n$, popularly referred to as L-statistics. L-statistics are usually represented as
\begin{equation}\label{lstat}
L_n =\displaystyle \sum_{i=1}^n a_{i,n}W_{(i)}, \quad 1 \leq i \leq n,
\end{equation}
where $a_{i,n}$ is a triangular array of constants, referred to as \emph{weights}. A wide variety of limiting results on L-statistics have been derived over the years. We direct the interested reader to \cite{ABN} for a good source of results and relevant references. The asymptotic properties of these objects have been determined under suitable regularity conditions, albeit usually not too stringent, nevertheless disconcerting on occasions in practice. In this paper, we examine one such occasion, wherein we are faced with an L-statistic---the ECF---whose weights are not sufficiently smooth. As a consequence, asymptotic normality and a functional limit theorem do not follow readily.

\cite{JH2}, in his elegant paper derived asymptotic distributions of clustering criteria. He employed, what he referred to as the \emph{split function}, in deriving the limiting results.  The ECF, $G_n$, arises in a natural manner as the empirical counterpart of a certain functional of his split function when we are concerned with random variables having common invertible distribution function. The ECF is an interesting probabilistic object in its own right and being a linear function of the data, offers an advantage over Hartigan's quadratic criterion function in terms of being amenable to extension to more interesting settings---namely clustering in higher dimensions and clustering of dependent observations.

The properties of the $k$-means clustering procedure for the univariate and the multivariate cases have been investigated extensively. \cite{DP1}, \cite{DP2} proved strong consistency and asymptotic normality results in the univariate case. \cite{SGJ} proved some weak limit theorems under non regular conditions for the univariate case. With the intention of having a more robust procedure for clustering, \cite{garcia1}, \cite{garcia2} propose the trimmed $k$-means clustering and provide a central limit theorem for the multivariate case.  In this paper, we prove consistency, a central limit theorem and also an invariance principle for our criterion function $G_n$, which is not in a form amenable for the usual representation of an L-statistic; nor, are its weights sufficiently smooth for the applicability of existing results.
\section{Empirical Cross-over Function} \label{crossover}
In this section, we introduce the necessary constructs from clustering techniques from which we develop the ECF. Let $W_1,W_2,\dots,W_n$  be i.i.d. random variables with continuous cumulative distribution function $F$. We make the following assumptions.
\begin{description}
\item [$A1$.] $F$ is invertible for $0<p<1$ and absolutely continuous with density $f$.
\item [$A2$.] $E(W_1^2)<\infty$.
\item [$A3$.] For $0<p<1$, $F$ is twice differentiable at $F^{-1}(p)$.
\end{description}
It is fairly common to encounter invertible distribution functions in applications. For example, models  in finance possess strictly increasing distribution functions usually guaranteed by the additive ``sort of Gaussian" noise from the Ito integral component which smooths and removes both jumps and flat areas of the distribution function.  \\

For $0<p<1$, consider the the \emph{split function} of $F^{-1}$ at $p$, as defined in \cite{JH2},
\begin{eqnarray*}
 B(F^{-1},p) = p\mu_l^2+(1-p)\mu_u^2-\left(\int_0^1F^{-1}(q)dq\right)^2,
\end{eqnarray*}
where
\begin{align*}
\mu_l&=\frac{1}{p}\int_{q\leq p}F^{-1}(q)dq=\frac{1}{p}\displaystyle\int_{-\infty}^{F^{-1}(p)}wdF, \\
\mu_u&=\frac{1}{1-p}\int_{q> p}F^{-1}(q)dq=\frac{1}{1-p}\displaystyle\int_{F^{-1}(p)}^{\infty}wdF.
\end{align*}
One way to think of $B(F^{-1},p)$ is, as the `between cluster sum of squares', in the case where we are concerned with two clusters in one dimension. Therefore, the value of $p \in (0,1)$ maximizing this function, would determine the location at which data is split into two clusters. Let us denote that value as $p_0$ and $p_0$ is referred to as the \emph{split point} in \cite{JH2}. As pointed out in \cite{JH2}, the conditions that guarantee the existence and uniqueness of the split point, are unclear. Determination of the requisite conditions, alone, is worthy of further investigation. However, for the purposes of this paper, those conditions and the split point itself are not important. When $F$ is invertible, it is known that the split point $p_0$ solves
\begin{equation}
\label{splitone}
(\mu_u-\mu_l)[\mu_u+\mu_l-2F^{-1}(p)]=0,
\end{equation}
where the left side is the derivative of $B(F^{-1},p)$. 
Owing to the fact that $(\mu_u-\mu_l) > 0$ for all $0<p<1$, we are interested only in the zero of
\begin{equation}
\label{split}
G(p)= \mu_l+\mu_u-2F^{-1}(p),
\end{equation}
which we refer to as the \emph{cross-over function}. The empirical version of the cross-over function represents the primary object of this paper. At this juncture, for better exposition, we recall the definition of the ECF; for $0<p<1$, we have
\begin{equation}\label{G_function}
 G_n(p)=\frac{1}{k}\sum_{j=1}^k W_{(j)}-W_{(k)} +\frac{1}{n-k}\sum_{j=k+1}^n W_{(j)}-W_{(k+1)} \qquad \textrm{for  }\frac{k-1}{n}\leq p <\frac{k}{n}.
\end{equation}
\begin{remark}
Intuition about the ECF is useful here. The term `cross-over' arises owing to the observation that
\begin{eqnarray*}
 G_n\left(\frac{1}{n}\right) &=& W_{(1)}-W_{(1)}+\frac{1}{n-1}\sum_{j=2}^n W_{(j)}-W_{(2)}\quad \geq 0,\\
  G_n\left(\frac{n-1}{n}\right)&=& \frac{1}{n-1}\sum_{j=1}^{n-1} W_{(j)}-W_{(n-1)}+W_{(n)}-W_{(n)} \quad \leq 0,
\end{eqnarray*}
and the function crosses over $0$ at some $1<k<n-1$. If $k^*$ is the index at which $G_n$ crosses over, then $W_{(k^*)}$ represents the datum at which the data is split leading to the formation of two clusters.  The term, $\frac{1}{k}\sum_{j=1}^k W_{(j)}-W_{(k)}$, can be thought of as a `distance' between the mean of the first $k$ observations, arranged in increasing order, and their maximum value; the term, $\frac{1}{n-k}\sum_{j=k+1}^n W_{(j)}-W_{(k+1)}$, represents the `distance' between the mean of the last $k$ observations and their minimum.
\end{remark}
\begin{remark}
The function $G_n$ is a linear combination of order statistics $W_{(i)}, 1 \leq i \leq n$ and hence an L-statistic. In the representation of an L-statistic $L_n$ shown in (\ref{lstat}), if the weights $a_{i,n},1 \leq i \leq n$, are of the form $\frac{1}{n}J\left(\frac{i}{n+1}\right)$, where $J(u),0 <u<1$, is the \emph{weight function}, then it is possible to obtain an equivalent representation as
\begin{equation*}
L_n=\frac{1}{n}\displaystyle \sum_{i=1}^n J\left(\frac{i}{n+1}\right)W_{(i)}.
\end{equation*}
The form of the weights $a_{i,n}$ represent the smoothness condition which guarantees asymptotic normality (See for e.g., \cite{ABN}, page 227 or \cite{van}, page 318). Unfortunately, $G_n$ cannot be represented in this form, since it has `bad' weights, in the following sense: for $0<p<1$, we see that the order statistics $W_{(\round)}$ and $W_{(\round+1)}$ have weights $\frac{1}{\round}-1$ and $\frac{1}{\roundp}-1$, respectively. This clearly violates the smoothness condition rendering the usage of existing results inappropriate.
\end{remark}
\begin{remark}
Observe that for a fixed $0<p<1$,
\begin{align*}
\frac{1}{k}\sum_{j=1}^k W_{(j)}-W_{(k)} &= \frac{1}{\round}\sum_{j=1}^{\round} W_{(j)}-W_{(\round)}, \\
\frac{1}{n-k}\sum_{j=k+1}^n W_{(j)}-W_{(k+1)} &= \frac{1}{\roundp}\sum_{j=\round+1}^{n} W_{(j)}-W_{(\round +1)},
\end{align*}
where $\lceil x\rceil$ represents the smallest integer not less than $x$. For a fixed $p \in (0,1)$, the sums shown above are \emph{trimmed} sums. More precisely, since $\frac{\round}{n} \rightarrow p$ and $\frac{\roundp}{n} \rightarrow 1-p$, they represent the case of \emph{heavy} trimming, asymptotics for which are well known (see for e.g., \cite{maller} and \cite{stigler}). Unfortunately, the two order statistics, $W_{(\round)}$ and $W_{(\roundp)}$, represent a formidable obstacle in the use of existing results regarding asymptotic normality of heavily trimmed sums. The function $G_n$ is hence some sort of a combination of heavily trimmed sums and intermediate order statistics, and asymptotic results for such a combination, to our knowledge, are yet to developed.
\end{remark}
\section{Limit thereoms for $G_n$}\label{results}
In this section, we prove the main results on the asymptotic behavior of the sample cross-over function $G_n$.
\begin{theorem}\label{th:consistency}
Under the assumptions $A1$ and $A2$ as $n\rightarrow \infty$,
\[
G_n(p)\overset{P}\rightarrow G(p).
\]
\end{theorem}
\begin{proof}
{Because we only need to prove consistency for individual components of the ECF, it is a relatively easy exercise. However, for a purpose of completeness and in order to introduce notation and ideas that will be used in the proof of the subsequent theorem, we decided to provide a detailed proof of the law of large numbers for $G_n$.}

For $0<p<1$, it is well known that $W_{(\lceil np \rceil)} \overset{P}\rightarrow F^{-1}(p)$ at points of continuity of $F^{-1}$. It is also the case that $W_{(\round+1)} \overset{P}\rightarrow F^{-1}(p)$, since the condition for $k_n$-th order statistic $W_{(k_n)}$ to be consistent for $F^{-1}(p)$  is that $\frac{k_n}{n}\rightarrow p$ (see for instance, \cite{van}). Let us define
\begin{equation*}
r_n=\frac{1}{n}\displaystyle \sum_{i=1}^n \less,
\end{equation*}
 where $\mathcal{I}_A$ is the indicator function of the set $A$. By the strong law of large numbers, $r_n \rightarrow p$ w.p.1. Now,
\begin{equation*}
  \frac{1}{k}\displaystyle \sum_{i=1}^{k}W_{(i)}= \frac{1}{\lceil np \rceil}\displaystyle \sum_{i=1}^{\lceil np \rceil}W_{(i)}.
\end{equation*}
Therefore,
\begin{equation*}
\frac{1}{k}\displaystyle \sum_{i=1}^{k}W_{(i)}= \frac{1}{\lceil np \rceil}\left[\displaystyle \sum_{i=1}^{\lceil nr_n \rceil}W_{(i)} + \displaystyle \sum_{i=\lceil nr_n \rceil+1}^{\lceil np \rceil}W_{(i)}\right]  .
\end{equation*}
It is clear here that if $\lceil nr_n \rceil+1 >\round$, 
the upper and lower limits of the second sum are interchanged with a negative sign.

The random sum
\begin{align*}
\frac{1}{\lceil np \rceil}\left|\displaystyle \sum_{i=\lceil nr_n \rceil+1}^{\lceil np \rceil}W_{(i)}\right| &\leq \frac{1}{\lceil np \rceil}\displaystyle \sum_{i=\lceil nr_n \rceil+1}^{\lceil np \rceil}|W_{(i)}|\\
&\leq  \frac{1}{\round}|\round-\rround|\left(|W_{(\round)}|+|W_{(\rround)}|\right).
\end{align*}
Recall that $r_n=p+O_p(n^{-1/2})$ and hence $|W_{\round}|$ and $|W_{\rround}|$ converge in probability to $F^{-1}(p)$(see \cite{van}, page 308) and $|r_n-p|\overset{P}\rightarrow 0$. Consequently, we have that
\begin{equation*}
\frac{1}{\round}\displaystyle \sum_{i=\rround+1}^{\round}W_{(i)}\overset{P}\rightarrow 0.
\end{equation*}

However,
\begin{equation*}
 \frac{1}{\round}\displaystyle \sum_{i=1}^{\rround}W_{(i)}= \frac{1}{\round}\displaystyle \sum_{i=1}^nW_i\less \overset{P}{\rightarrow} \mu_l \quad,
\end{equation*}
by the law of large numbers of i.i.d random variables. As a consequence, $ \frac{1}{k}\displaystyle \sum_{i=1}^{k}W_{(i)} \overset{P}{\rightarrow}\mu_l$. In similar fashion we note that $  \frac{1}{n-k}\displaystyle \sum_{i=k+1}^{n}W_{(i)} \overset{P}{\rightarrow} \mu_u$ for all $0<p<1$ and $\frac{k-1}{n}\leq p < \frac{k}{n}$. Combining the above two convergences with the convergence of $W_{(\round)}$ and $W_{(\round+1)}$ to their identical limits, we have that, $G_n(p)\overset{P}\rightarrow G(p)$ for each $0<p<1$.
\end{proof}
\begin{remark}
It is worthwhile to note that the trimmed (at the random level)  sum $ \sum_{i=1}^{\rround}W_{(i)}$ is exactly equal to the truncated sum $\sum_{i=1}^nW_i\less$, which is the sum of i.i.d. random random variables. This subtle relationship is greatly convenient in our proofs.
\end{remark}

For ease of notation, let us define for $0<p<1$,
\begin{align*}
\theta_p=&\frac{1}{p}W_1\lessone-\frac{1}{p}F^{-1}(p)\lessone\\
&+\frac{1}{1-p}W_1\moreone-\frac{1}{1-p}F^{-1}(p)\moreone\\
&+\frac{2\lessone}{f(F^{-1}(p))},
\end{align*}
and $U_n(p)=\sqrt{n}(G_n(p)-G(p))$ for $0<a\leq p \leq b<1$.
\begin{theorem}
Under assumptions $A1-A3$ as $n\rightarrow \infty$,
\[
\sqrt{n}\left(G_n(p)-G(p)\right)\overset{d}\rightarrow N\left(0,\sigma\right),
\]
where $\sigma=Var(\theta_{p})$. Furthermore,
\begin{equation*}
U_n(p)\Rightarrow U(p),
\end{equation*}
in the Skorohod space $D[a,b]$, $0<a<b<1$ equipped with the $J_1$ topology, where $U(p)$ is a Gaussian process with mean $0$ and covariance given by
\[
Cov(U(p),U(q))=Cov(\theta_p,\theta_q).
\]
\end{theorem}
\begin{proof}
The trick used in proving the asymptotic normality of $G_n$ is to consider mean-zero asymptotics of its individual components and by the use of \emph{Bahadur's representation} for sample quantiles, rewrite $G_n$ as a sum of i.i.d random variables and an error term, which goes to zero at an appropriate rate. This would then pave the way for the usage of standard results.

More specifically, first note that for $0<p<1$, and each $i=1,\dots,n$, $$E\left(W_i \less\right)=p\mu_l$$ and $$E\left(W_i \more\right)=(1-p)\mu_u.$$
Observe that, for $\frac{k-1}{n} \leq p <\frac{k}{n}$,
\begin{align*}
\sqrt{n}\left[\frac{1}{k}\displaystyle \sum_{i=1}^{k}W_{(i)}-np\mu_l\right] &= \frac{\sqrt{n}}{\round}\left[\displaystyle \sum_{i=1}^{\round}W_{(i)}-np\mu_l\right]\\
&= \frac{\sqrt{n}}{\round}\left[\displaystyle \sum_{i=1}^{\rround}W_{(i)}+\displaystyle \sum_{i=\rround+1}^{\round}W_{(i)}-np\mu_l\right]\\
&=\frac{\sqrt{n}}{\round}\left[\displaystyle \sum_{i=1}^{\rround}W_{(i)}-np\mu_l\right]+\frac{\sqrt{n}}{\round}\left[\displaystyle \sum_{i=\rround+1}^{\round}\left(W_{(i)}-F^{-1}(p)\right)\right]\\
&\quad+\frac{\sqrt{n}}{\round}F^{-1}(p)\left(\round-\rround\right).
\end{align*}
Now note that
\begin{align*}
\frac{\sqrt{n}}{\round}&\left|\displaystyle \sum_{i=\rround+1}^{\round}W_{(i)}-F^{-1}(p)\right| \\
&\qquad \qquad\leq
\frac{\sqrt{n}}{\round}|\round-\rround|\textrm{ max}\left(|W_{(\round)}-F^{-1}(p)|,|W_{(\rround)}-F^{-1}(p)|\right).
\end{align*}
By the same argument used in the proof of Theorem \ref{th:consistency}, $|W_{(\round)}-F^{-1}(p)|\overset{P}\rightarrow 0$ and $|W_{(\rround)}-F^{-1}(p)|\overset{P}\rightarrow 0$. By the central limit theorem for i.i.d random variables, $\sqrt{n}|p-r_n|$ is asymptotically normal and hence bounded in probability. Consequently,
\begin{equation*}
\frac{\sqrt{n}}{\round}\left|\displaystyle \sum_{i=\rround+1}^{\round}W_{(i)}-F^{-1}(p)\right|
\overset{P} \rightarrow 0 \quad .
\end{equation*}
Next, recall that
\begin{align*}
\displaystyle \sum_{i=1}^{\rround}W_{(i)}  &=   \displaystyle \sum_{i=1}^n W_i\less,\\
nr_n&=\displaystyle \sum_{i=1}^n\less.
\end{align*}
 Therefore,
\begin{align*}
\sqrt{n}\left[\frac{1}{\round}\displaystyle \sum_{i=1}^{\round}W_{(i)}-np\mu_l\right] &= \frac{1}{p}
\Big[\frac{1}{\sqrt{n}}\displaystyle \sum_{i=1}^n (W_i \less -p\mu_l)\\
        & \qquad+ \frac{1}{\sqrt{n}}F^{-1}(p)\displaystyle \sum_{i=1}^n(p-\less)\Big]\\
        & \qquad+o_{p}(1)\\
        &= \sqrt{n}\bar{\xi}+o_{p}(1),
\end{align*}
where
\begin{equation*}
\xi_i=\frac{1}{p}\left[W_i\less-F^{-1}(p)\less -(p\mu_l-pF^{-1}(p))\right]
\end{equation*}
are i.i.d random variables for $i=1,\dots,n$ and $\bar{\xi}=\frac{1}{n}\displaystyle \sum_{i=1}^n \xi_i$.

Using a similar argument, we can claim that
\begin{equation*}
\sqrt{n}\left[\frac{1}{\roundp}\displaystyle \sum_{i=\round+1}^{n}W_{(i)}-n(1-p)\mu_u\right]=\sqrt{n}\bar{\tau}+o_{p}(1),
\end{equation*}
where
\begin{equation*}
\tau_i=\frac{1}{1-p}\left[W_i\more-F^{-1}(p)\more -((1-p)\mu_u-(1-p)F^{-1}(p))\right]
\end{equation*}
are i.i.d random variables and $\bar{\tau}=\frac{1}{n}\displaystyle \sum_{i=1}^n \tau_i$. That takes care of the two trimmed sums.

Next, we  turn our attention to the two quantiles $W_{(k)}$ and $W_{(k+1)}$ or equivalently $W_{(\round)}$ and $W_{(\round+1)}$ for $\frac{k-1}{n} \leq p <\frac{k}{n}$. Using the Bahadur representation for sample quantiles (see \cite{bahadur}), justified by assumptions $A1$ and $A3$ , we have
\begin{equation*}
\sqrt{n}\left(W_{(\round)}-F^{-1}(p)\right)=\sqrt{n}\left(W_{(\round+1)}-F^{-1}(p)\right)=\sqrt{n}\bar{\kappa}+o_p(1),
\end{equation*}
where $\bar{\kappa}=\frac{1}{n}\displaystyle \sum_{i=1}^n\kappa_i$, and
\begin{equation*}
\kappa_i=\frac{p-\less}{f(F^{-1}(p))}
\end{equation*}
are i.i.d random variables.

We are now in a situation where  for $0<p<1$, $\sqrt{n}(G_n(p)-G(p))$ has been expressed as sums of i.i.d random variables along with an error term which is $o_p(1)$. That is,
\begin{equation*}
\sqrt{n}\left(G_n(p)-G(p)\right)  = \displaystyle \sum_{i=1}^n \frac{Z_i}{\sqrt{n}}+o_p(1) ,
\end{equation*}
 where $Z_i=\xi_i+\tau_i-2\kappa_i$ are i.i.d random variables. The advantage of this representation lies in the fact that we are now allowed to examine $G_n$ without having to concern ourselves with the correlations between its individual components. The  representation ensures that the effect of the correlations is of order as that of the error term or smaller and can hence be safely disregarded.
 Consequently, by the central limit theorem for i.i.d random variables
 \begin{equation*}
 \sqrt{n}\left(G_n(p)-G(p)\right) \overset{d} \rightarrow N\left(0,\sigma\right).
 \end{equation*}

We can now turn our attention to the functional limit of the process $U_n$. Since we are interested in the behavior of $G_n$ for $0<p<1$ and in particular the point at which it crosses zero, we restrict ourselves to examining the behavior of $U_n$ in the closed interval $[a,b]$ where $a$ and $b$ are constants bounded away from $0$ and $1$ respectively. Notice that $U_n$ is a natural random element of the Skorohod space $D[a,b]$. It is straightforward to note that by virtue of our representation of $\sqrt{n}(G_n(p)-G(p))$, for each $p$, as {\it a sum of i.i.d. random variables plus an error term of order $o_p(1)$}, we can employ the central limit theorem for random vectors and obtain
\[
\sqrt{n}\left(U_n(p_1)-U(p_1),\dots,U_n(p_k)-U(p_k)\right)\overset{d}\rightarrow N\left(\mathbf{0},\Sigma \right),
\]
where $k$ is a finite positive integer and for $i,j=1,\dots,k$,
$\Sigma=\Bigl(\sigma_{ij}\Bigr)$ with
\[
  \sigma_{ij} = \left\{
  \begin{array}{l l}
    Var(\theta_{p_i}) & \quad \text{if $i=j$}\\
    Cov(\theta_{p_i},\theta_{p_j}) & \quad \text{if $i \neq j$}.\\
  \end{array} \right.
\]
Now, if we can show that the sequence $U_n$ is tight, we then have the required convergence to $U$ (see \cite{bill}, for the necessary arguments). We set about proving tightness in an indirect way, as opposed to the usual method of showing that $U_n$ concentrates on a compact set in $D[a,b]$ with high probability. Consider the components of $U_n$
\begin{align*}
&U^n_1=\sqrt{n}\left(\frac{1}{\round}\bigsum_{i=1}^{\round}\W-\frac{1}{p}\displaystyle \int_0^p F^{-1}(q)dq \right),\\
&U^n_2=\sqrt{n}\left(W_{(\round)}-F^{-1}(p)\right),\\
&U^n_3=\sqrt{n}\left(\frac{1}{\roundp}\bigsum_{i=\round+1}^{n}\W-\frac{1}{1-p}\displaystyle \int_p^1 F^{-1}(q)dq \right),\\
&U^n_4=\sqrt{n}\left(W_{(\round+1)}-F^{-1}(p)\right).
\end{align*}
It is interesting that for every $U_i^n$, $i=1,\dots,4$ the functional CLT is an established result. However, the weak convergence of the individual components $U_i^n$ does not automatically guarantee weak convergence for the sum of the components. But at this point  we need only the tightness.  Since the sum of compact sets is a compact set again, it is easy to show that if each component is tight then it is indeed true that the sum is tight with respect to the Skorohod metric on $D[a,b]$.  Now, note that $U_2^n$ and $U_4^n$ are quantile processes and converge weakly to a Gaussian process (see p. 308, \cite{van}) in $D[a,b]$. Using the result from \cite{KM}, we can claim that $U_1^n$ and $U_3^n$ also converge weakly to a limit process in $D[a,b]$. This proves that $U_i^n$ is relatively compact for each $i$. Now, since $D[a,b]$ is complete and separable with respect to the Skorohod metric (see p.115, \cite{bill}), using the converse of Prohorov's theorem (see p.37 \cite{bill}) we can claim that each $U_i^n$ for $i=1,\dots,4$  is tight and, therefore,  $U_n=U_1^n+U_2^n+U_3^n+U_4^n$  is tight in $D[a,b]$ equipped with the $J_1$-topology.

 \end{proof}

 We now provide verification of our asymptotic results regarding consistency and normality by considering two examples. In both the examples we first generate $1000$ random variables $T_n=\sqrt{n}\left(G_n(0.5)-G(0.5)\right)$ and obtain the simulated mean and the variance. In order to verify asymptotic normality, we generate again $100$ random variables $T_n$. This is done for different samples sizes $n$ and results are tabulated.
 \begin{example}\label{example1}
If $W_1,W_2,\dots,W_n$ are i.i.d. $N(0,1)$, then it can be ascertained quite easily that $G(0.5)=0$ and $\sigma=2\pi-4\approx2.2831$. The numbers tabulated below offer satisfactory evidence about the accuracy of our results.
\end{example}
\begin{example}
In this example, we consider $W_1,W_2,\dots,W_n$ to be i.i.d. exponential random variables with mean $1$. This represents the archetypal case of a skewed distribution and we again check for the accuracy of our results. In this case, $G(0.5)=2(1-\ln2)\approx0.6137$ and $\sigma=8(1-\ln2)\approx2.4548$. The numbers in the tables below provide further corroborative evidence for our limiting results.
\end{example}

\begin{table*}[h]
\begin{center}
\caption{Simulated means and variances for different sample sizes.}
\vspace{2mm}
\begin{tabular}{|c||c|c|c|c|c|c|}
\hline
\textbf{Random variables}&\multicolumn{3}{|c|}{$\mathbf{N(0,1)}$}&\multicolumn{3}{|c|}{$\mathbf{Exp(1)}$}\\
\hline
Sample sizes& $n=100$ & $n=1000$& $n=10000$& $n=100$ & $n=1000$& $n=10000$ \\ \hline
Simulated Mean &$-0.017$& $0.018$&$0.002$ &$-0.041$& $-0.014$&$0.0019$ \\ \hline
Simulated Variance &$2.407$& $2.324$&$2.296$ &$2.491$& $2.463$&$2.452$ \\ \hline
\end{tabular}
\end{center}
\end{table*}

\begin{table}[h]
\begin{center}
\caption{p-values for Kolmogorov-Smirnov test for normality}
\vspace{2mm}
\begin{tabular}{|c||c|c|}
\hline
\textbf{Random variables}& {$\mathbf{N(0,1)}$} & {$\mathbf{Exp(1)}$}\\ \hline
$n=100$& $0.751$ & $0.8786$ \\ \hline
$n=1000$& $0.12$ & $0.2174$ \\ \hline	
$n=10000$& $0.391$ &$0.9955$ \\	\hline

\end{tabular}
\end{center}
\end{table}

 \section{Concluding Remarks}\label{conclusion}

 Despite being an L-statistic, the asymptotic properties of the ECF cannot be studied using existing machinery owing to the fact that its weights are not smooth.  Asymptotic results for heavily trimmed sums are inapplicable to our problem due the presence of the two order statistics, $W_{(k)}$ and $W_{(k+1)}$, with unfriendly weights. The centered ECF,  however, can be expressed as a sum of i.i.d. random variables and an error term, which goes to zero at an appropriate rate, by the use of a subtle trick involving truncated sums and the Bahadur representation for sample quantiles. Owing to this, the CLT follows immediately and what remains is to show that the centered process satisfies the tightness condition for the functional CLT.

Note that the ECF is invariant with respect to shift in the distribution of $W$'s, but is linear with respect to scaling. If we introduce statistic $p_n$ (the empirical split point) that `solves', in some appropriate sense, the  equation $$G_n(p)=0,$$ then this statistic is invariant with respect to both shifting and scaling (as it should be, because the clustering problem is invariant with respect to linear transformations), and potentially can be used to design a clustering test.

The asymptotics of $p_n$ is the next natural question, which is the focus of the work in \cite{KB2}. The Central Limit Theorem for $G_n$, proved in this paper, constitutes a very important step towards the solution of determining the asymptotic behavior of $p_n$. According to a general plan outlined in \cite{serfling} p. 95, we can conjecture that $$p_n\approx p_0-G_n(p_0)/G'(p_0),$$ where $p_0$ is a theoretical split point.  However, the rigorous proof of this statement requires significant efforts.
\bibliography{ref}
\bibliographystyle{plainnat}

\end{document}